\documentclass{elsart}
\usepackage{amssymb,amsfonts,amsbsy}

\newcommand{\K}{\mathbb{K}}
\newcommand{\F}{\mathbb{F}}
\newcommand{\Fix}{\mathrm{Fix}}

\begin{document}

\begin{frontmatter}
\title{Building counterexamples to  generalizations for rational functions of Ritt's decomposition Theorem}
\author[unican]{Jaime Gutierrez\thanksref{x}}
\author[concordia]{David Sevilla\thanksref{x}}
\address[unican]{Dpto. de Matem\'aticas, Estad\'{\i}stica y Computaci\'on, Universidad de Cantabria, E--39071 Santander, Spain}
\address[concordia]{Dpt. of Computer Science and Software Engineering, University of Concordia, Montreal, Canada}
\thanks[x]{Partially supported by Spain Ministry of Science grant MTM2004-07086}

\begin{abstract}
The classical Ritt's Theorems state several properties of univariate polynomial decomposition. In
this paper we present new counterexamples to the first Ritt theorem, which states the equality of
length of decomposition chains of a polynomial, in the case of rational functions. Namely, we
provide an explicit example of a rational function with coefficients in $\Qset$ and two
decompositions of different length.

Another aspect is the use of some techniques that could allow for other counterexamples, namely,
relating groups and decompositions and using the fact that the alternating group $A_4$ has two
subgroup chains of different lengths; and we provide more information about the generalizations of
another property of polynomial decomposition: the stability of the base field. We also present an
algorithm for computing the fixing group of a rational function providing the complexity over the
rational number field.
 \end{abstract}
\end{frontmatter}

\section{Introduction}

The starting point is the decomposition of polynomials and rational functions in one variable.
First we will define the basic concepts of this topic.

\begin{defn}
If $f=g\circ h$, $f,g,h\in\K(x)$, we call this a \emph{decomposition of $f$ in $\K(x)$} and say
that $g$ is a \emph{component on the left} of $f$ and $h$ is a \emph{component on the right} of
$f$. We call a decomposition \emph{trivial} if any of the components is a unit with respect to
decomposition.

Given two decompositions $f=g_1\circ h_1=g_2\circ h_2$ of a rational function, we call them
\emph{equivalent} if there exists a unit $u$ such that
\[h_1=u\circ h_2\,,\qquad g_1=g_2\circ u^{-1},\]
where the inverse is taken with respect to composition.

Given a non--constant $f$, we say that it is \emph{indecomposable} if it is not a unit and all its
decompositions are trivial.

We define a \emph{complete} decomposition of $f$ to be $f=g_1\circ\cdots\circ g_r$ where $g_i$ is
indecomposable. The notion of equivalent complete decompositions is straightforward from the
previous concepts.

Given a non--constant rational function $f(x)\in\K(x)$ where $f(x)=f_N(x)/f_D(x)$ with
$f_N,f_D\in\K[x]$ and $(f_N,f_D)=1$, we define the \emph{degree} of $f$ as
\[\deg\,f=\max\{\deg\,f_N,\ \deg\,f_D\}.\]

We also define $\deg\,a=0$ for each $a\in\K$.
\end{defn}

\begin{rem}
From now on, we will use the previous notation when we refer to the numerator and denominator of a
rational function. Unless explicitly stated, we will take the numerator to be monic, even though
multiplication by constants will not be relevant.
\end{rem}

The first of Ritt's Theorems states that all the decomposition chains of a polynomial that
satisfies a certain condition have the same length.  Here we explore new techniques related to this,
and include a counterexample in $\Qset(x)$.

Another result in this fashion states that if a polynomial is indecomposable in a certain
coefficient field, then it is also indecomposable in any extension of that field. This is also
false for rational functions, see \cite{DW74} and \cite{AGR95}. We look for bounds for the degree
of the extension in which we need to take the coefficients if a rational function with
coefficients in $\Qset$ has a decomposition in a larger field. In this paper we present a
computational approach to this question and our conclusions.

In Section 2 we study how to compute bounds for the minimal field that contains all the
decompositions of a given rational function. In Section 3 we introduce several definitions and
properties of groups related to rational functions, which we use in Section 4 to discuss the
number of components in the rational case. In particular, we present an algorithm for computing
fixing group of a rational function and we provide the complexity over the rational number field.
Finally, in Section 4 we present an example of a degree 12 rational function with coefficients in
$\Qset$ and two decompositions of different length; as far as we know this is the first example in
$\Qset$ of this kind.

\section{Extension of the coefficient field}

Several algorithms for decomposing univariate rational functions are known, see for instance
\cite{Zip91} and \cite{AGR95}. In all cases, the complexity of the algorithm grows enormously when
the coefficient field is extended. A natural question about decomposition is whether it depends on
the coefficient field, that is, the existence of polynomials or rational functions that are
indecomposable in $\K(x)$ but have a decomposition in $\F(x)$ for some extension $\F$ of $\K$.
Polynomials behave well under certain conditions, however in the rational case this is not true.
We will try to shed some light on the rational case.

\begin{defn}
$f\in\K[x]$ is \emph{tame} when $\mathrm{char}\ \K$ does not divide $\deg\,f$.
\end{defn}

The next theorem shows that tame polynomials behave well under extension of the coefficient field,
see \cite{Gut91}. It is based on the concept of \emph{approximate root} of a polynomial, which
always exists for tame polynomials, and is also the key to some other structural results in the
tame polynomial case.

\begin{thm}\label{indesc-ext}
Let $f\in\K[x]$ be tame and $\F\supseteq\K$. Then $f$ is indecomposable in $\K[x]$ if and only if
it is indecomposable in $\F[x]$.
\end{thm}

The next example, presented in \cite{AGR95}, shows that the previous result is false for rational
functions.

\begin{exmp}\label{contraej-extension-rac}
Let
\[f=\frac{\omega^3x^4-\omega^3x^3-8x-1}{2\omega^3x^4+\omega^3x^3-16x+1}\]
where $\omega\not\in\Qset$ but $\omega^3\in\Qset\setminus\{1\}$. It is easy to check that $f$ is
indecomposable in $\Qset(x)$. However, $f=f_1\circ f_2$ where
\[f_1=\frac{x^2+(4-\omega)x-\omega}{2x^2+(8+\omega)x+\omega}\,,\qquad f_2=\frac{x\omega(x\omega-2)}{x\omega+1}\,.\]
\end{exmp}

We can pose the following general problem:

\begin{prob}
Given a function $f\in\K(x)$, compute a minimal field $\F$ such that every decomposition of $f$
over an extension of $\K$ is equivalent to a decomposition over $\F$.
\end{prob}

It is clear that, by composing with units in $\F(x)\supseteq\K(x)$, we can always turn a given
decomposition in $\K(x)$ into one in $\F(x)$. Our goal is to minimize this, that is, to determine
fields that contain the smallest equivalent decompositions in the sense of having the smallest
possible extension over $\K$.

Given a decomposition $f=g(h)$ of a rational function in $\K(x)$, we can write a polynomial system
of equations in the coefficients of $f$, $g$ and $h$ by equating to zero the numerator of
$f-g(h)$. The system is linear in the coefficients of $g$. Therefore, all the coefficients of $g$
and $h$ lie in some algebraic extension of $\K$. Our goal is to find bounds for the degree of the
extension $[\F:\K]$ where $\F$ contains, in the sense explained above, all the decompositions of
$f$.

One way to find a bound is by means of a result that relates decomposition and factorization. We
state the main definition and theorems here, see \cite{GRS01} for proofs and other details.

\begin{defn}
A rational function $f\in\K(x)$ is in \emph{normal form} if $\deg\,f_N>\deg\,f_D$ and $f_N(0)=0$
(thus $f_D(0)\neq 0$).
\end{defn}

\begin{thm}\label{exist-unids}$ $

\item[(i)] Given $f\in\K(x)$, if $\deg\,f<|\K|$ then there exist units $u,v$ such that $u\circ f\circ
v$ is in normal form.

\item[(ii)] If $f\in\K(x)$ is in normal form, every decomposition of $f$ is equivalent to one where both
components are in normal form.
\end{thm}

We will analyze the complexity of finding the units $u$ and $v$ later.

\begin{thm}\label{divide-normal-univ}
Let $f=g(h)$ with $f,g,h$ in normal form. Then $h_N$ divides $f_N$ and $h_D$ divides $f_D$.
\end{thm}

This result provides the following bound.

\begin{thm}\label{cota-grado-ext}
Let $f\in\K(x)$ and $u_1,u_2$ be two units in $\K(x)$ such that
$g=u_1\circ f\circ u_2$ is in normal form. Let $\F$ be the splitting
field of $\{g_N,g_D\}$.  Then any decomposition of $f$ in $\K'(x)$,
for any  $\K'\supset\K$ is equivalent to a decomposition in $\F(x)$.
\end{thm}
\begin{pf}
By Theorems \ref{exist-unids} and \ref{divide-normal-univ}, every decomposition of $g$ is
equivalent to another one, $g=h_1\circ h_2$, where the numerator and denominator of $h_2$ divide
those of $g$, thus the coefficients of that component are in $\F$. As the coefficients of $h_1$
are the solution of a linear system of equations whose coefficients are polynomials in the
coefficients of $g$ and $h_2$, they are also in $\F$. We also have $u_1,u_2\in\K(x)$, therefore
the corresponding decomposition of $f$ lies in the same field. \qed
\end{pf}

This bound, despite being of some interest because its generality and simplicity, is far from
optimal. For example, for degree 4 we obtain $[\F:\K]\leq 3!\cdot 3!=36$.
The following theorem completes Example \ref{contraej-extension-rac}.

\begin{thm}
Let $f\in\Qset(x)$ of degree 4. If $f=g(h)$ with $g,h\in\overline{\Qset}(x)$, there exists a field
$\K$ with $\Qset\subset\K\subset\overline{\Qset}$ and a unit $u\in\K(x)$ such that
$g(u^{-1}),u(h)\in\K(x)$ and $[\K:\Qset]\leq 3$.
\end{thm}

The proof
is a straightforward application of  Gr\"obner bases and the
 well--known Extension Theorem, see for instance  \cite{CLO97}.
\section{Fixing group and fixed field}
In this section we introduce several simple notions from classical Galois theory. Let
$\Gamma(\K)=\mathrm{Aut}_\K\K(x)$ (we will write simply $\Gamma$ if there can be no confusion
about the field). The elements of $\Gamma(\K)$ can be identified with the images of $x$ under the
automorphisms, that is, with M\"obius transformations (non--constant rational functions of the
form $(ax+b)/(cx+d)$), which are also the units of $\K(x)$ under composition.

\begin{defn} $ $

\item[(i)] Let $f\in\K(x)$. We define $G(f)=\{u\in\Gamma(\K):\ f\circ u=f\}$.

\item[(ii)] Let $H<\Gamma(\K)$. We define $\Fix(H)=\{f\in\K(x):\ f\circ u=f\ \forall u\in H\}$.
\end{defn}

\begin{exmp}
\item[(i)] Let $f=x^2+\displaystyle\frac{1}{x^2}\in\K(x)$. Then $G(f)=\left\{x,-x,\displaystyle\frac{1}{x},-\displaystyle\frac{1}{x}\right\}$.

\item[(ii)] Let $H=\{x,ix,-x,-ix\}\subset\Gamma(\Cset)$. Then $\Fix(H)=\Cset(x^4)$.
\end{exmp}

These definitions correspond to the classical Galois correspondences (not bijective in general)
between the intermediate fields of an extension and the subgroups of its automorphism group, as
the following diagram shows:

\[\begin{array}{ccc}
    \K(x)  &  \longleftrightarrow  &  \{id\}  \\
    |  &  &  |  \\
    \K(f)  &  \longrightarrow  &  G(f)  \\
    |  &  &  |  \\
    \Fix(H)  &  \longleftarrow  &  H  \\
    |  &  &  |  \\
    \K  &  \longleftrightarrow  &  \Gamma  \\
\end{array}\]

\begin{rem}
As $\K(f)=\K(f')$ if and only if $f=u\circ f'$ for some unit $u$, we have that the application
$\K(f)\mapsto G(f)$ is well--defined.
\end{rem}

Next, we state several interesting properties of the fixed field and the fixing group.

\begin{thm}\label{H-inf-fin}
Let $H$ be a subgroup of $\Gamma$.

\item[(i)] $H$ is infinite $\Rightarrow \Fix(H)=\K$.

\item[(ii)] $H$ is finite $\Rightarrow\K\varsubsetneq\Fix(H)$, $\Fix(H)\subset\K(x)$ is a normal
extension, and in particular $\Fix(H)=\K(f)$ with $\deg\,f=|H|$.
\end{thm}

\begin{pf} $ $

(i) It is clear that no non--constant function can be fixed by infinitely many units, as these
must fix the roots of the numerator and denominator.

(ii) We will show constructively that there exists $f$ such that $\Fix(H)=\K(f)$ with
$\deg\,f=|H|$. Let $H=\{h_1=x,\ldots,h_m\}$. Let
\[P(T)=\prod_{i=1}^m (T-h_i)\ \in\ \K(x)[T].\]
We will see that $P(T)$ is the minimum polynomial of $x$ over $\Fix(H)\subset\K(x)$. A classical
proof of L\"uroth's Theorem (see for instance \cite{Wae64}) states that any non--constant
coefficient of the minimum polynomial generates $\Fix(H)$, and we are done.

It is obvious that $P(x)=0$, as $x$ is always in $H$. It is also clear that $P(T)\in\Fix(H)[T]$,
as its coefficients are the symmetric elementary polynomials in $h_1,\ldots,h_m$. The
irreducibility is equivalent to the transitivity of the action of the group on itself by
multiplication. \qed
\end{pf}

\begin{thm}\label{props-fix}$ $

\item[(i)] For any non--constant $f\in\K(x)$, $|G(f)|$ divides $\deg\,f$. Moreover, for any field $\K$
there is a function $f\in\K(x)$ such that $1<|G(f)|<\deg\,f$.

\item[(ii)] If $|G(f)|=\deg\,f$ then $\K(f)\subseteq\K(x)$ is normal. Moreover, if the extension
$\K(f)\subseteq\K(x)$ is separable, then
\[\K(f)\subseteq\K(x)\mbox{ is normal}\quad \Rightarrow\quad |G(f)|=\deg\,f.\]

\item[(iii)] Given a finite subgroup $H$ of $\Gamma$, there is a bijection between the subgroups of
$H$ and the fields between $\Fix(H)$ and $\K(x)$. Also, if $\Fix(H)=\K(f)$, there is a bijection
between the right components of $f$ (up to equivalence by units) and the subgroups of $H$.
\end{thm}

\begin{pf} $ $

(i) The field $\Fix(G(f))$ is between $\K(f)$ and $\K(x)$, therefore the degree of any generator,
which is the same as $|G(f)|$, divides $\deg\,f$. For the second part, take for example
$f=x^2\,(x-1)^2$, which gives $G(f)=\{x,1-x\}$ in any coefficient field.

(ii) The elements of $G(f)$ are the roots of the minimum polynomial of $x$ over $\K(f)$ that are
in $\K(x)$. If there are $\deg\,f$ different roots, as this number equals the degree of the
extension we conclude that it is normal.

If $\K(f)\subset\K(x)$ is separable, all the roots of the minimum polynomial of $x$ over $\K(f)$
are different, thus if the extension is normal there are as many roots as the degree of the
extension.

(iii) Due to Theorem \ref{H-inf-fin}, the extension $\Fix(H)\subset\K(x)$ is normal, and the
result is a consequence of the Fundamental Theorem of Galois.
\end{pf}

\begin{rem}
$\K(x)$ is Galois over $\K$ (that is, the only rational functions fixed by $\Gamma(\K)$ are the
constant ones) if and only if $\K$ is infinite. Indeed, if $\K$ is infinite, for each
non--constant function $f$ there exists a unit $x+b$ with $b\in\K$ which does not leave it fixed.
On the other hand, if $\K$ is finite then $\Gamma(\K)$ is finite too, an the proof of Theorem
\ref{H-inf-fin} provides a non--constant rational function that generates $\Fix(\Gamma(\K))$.
\end{rem}

Algorithms for computing several aspects of Galois theory can be found in \cite{Val95}.
Unfortunately, it is not true in general that $[\K(x):\K(f)]=|G(f)|$; there is no bijection
between intermediate fields and subgroups of the fixing group of a given function. Anyway, we can
obtain partial results on decomposability.

\begin{thm}\label{conds-indesc}
Let $f$ be indecomposable.
\item[(i)] If $\deg\,f$ is prime, then either $G(f)$ is cyclic of order $\deg\,f$, or it is trivial.

\item[(ii)] If $\deg\,f$ is composite, then $G(f)$ is trivial.
\end{thm}

\begin{pf} $ $

(i) If $1<|G(f)|<\deg\,f$, we have $\K(f)\subsetneq\K(\Fix(G(f)))\subsetneq\K(x)$ and any
generator of $\K(\Fix(G(f)))$ is a proper component of $f$ on the right. Therefore, $G(f)$ has
order either 1 or $\deg\,f$, and in the latter case, being prime, the group is cyclic.

(ii) Assume $G(f)$ is not trivial. If $|G(f)|<\deg\,f$, we have a contradiction as in (i). If
$|G(f)|=\deg\,f$, as it is a composite number, there exists $H\lneq G(f)$ not trivial, and again
any generator of $\Fix(H)$ is a proper component of $f$ on the right.
\end{pf}

\begin{cor}
If $f$ has composite degree and $G(f)$ is not trivial, $f$ is decomposable.
\end{cor}

Now we present algorithms to efficiently compute fixed fields and fixing groups.

The proof of Theorem \ref{H-inf-fin} provides an algorithm to compute a generator of $\Fix(H)$
from its elements.

\begin{alg}\label{alg-fixed-field} $ $
\begin{description}
    \item \emph{INPUT}: $H=\{h_1,\ldots,h_m\}<\Gamma(\K)$.

    \item \emph{OUTPUT}: $f\in\K(x)$ such that $\Fix(H)=\K(f)$.
\end{description}
\begin{description}
    \item \emph{A}. Let $i=1$.

    \item \emph{B}. Compute the $i$-th symmetric elementary function $\sigma_i(h_1,\ldots,h_m)$.

    \item \emph{C}. If $\sigma_i(h_1,\ldots,h_m)\not\in\K$, return $\sigma_i(h_1,\ldots,h_m)$. If
it is constant, increase $i$ and return to \textsc{B}.
\end{description}
\end{alg}

We illustrate this algorithm with the following example.

\begin{exmp}\label{ej-fixed-field2}
Let
\[H=\left\{\pm x\ ,\pm \frac{1}{x}\ ,\pm \frac{i(x+1)}{x-1}\ ,\pm \frac{i(x-1)}{x+1}\ ,
\pm \frac{x+i}{x-i}\ ,\pm \frac{x-i}{x+i}\right\}<\Gamma(\Cset).\]

Then
\[\begin{array}{rl}
  P(T)  =  &  T^{12}-\displaystyle{\frac{x^{12}-33x^8-33x^4+1}{x^2(x-1)^2(x+1)^2(x^4+2x^2+1)}}\ T^{10}-33\,T^8  \\
    &  +\ 2\ \displaystyle{\frac{x^{12}-33x^8-33x^4+1}{x^2(x-1)^2(x+1)^2(x^4+2x^2+1)}}\ T^6-33\,T^4  \\
    &  -\ \displaystyle{\frac{x^{12}-33x^8-33x^4+1}{x^2(x-1)^2(x+1)^2(x^4+2x^2+1)}}\ T^2+1\,.
\end{array}\]

Thus,
\[\Fix(H)=\Cset\left(\frac{x^{12}-33x^8-33x^4+1}{x^2(x-1)^2(x+1)^2(x^4+2x^2+1)}\right).\]

$H$ is isomorphic to $A_4$. It is known that $A_4$ has two complete subgroup chains of different
lengths:
\[\{id\}\subset C_2\subset V\subset A_4\,,\qquad \{id\}\subset C_3\subset A_4\,.\]

In our case,
\[\{x\}\subset\{\pm x\}\subset\left\{\pm x,\pm \frac{1}{x}\right\}\subset H\,,\qquad \{x\}\subset\left\{x,\frac{x+i}{x-i},\frac{i(x+1)}{x-1}\right\}\subset H\,.\]

Applying our algorithm again we obtain the following field chains:
\[\begin{array}{l}
  \Cset(f)\ \subset\ \Cset\left(x^2+\displaystyle\frac{1}{x^2}\right)\ \subset\ \Cset(x^2)\ \subset\ \Cset(x)\,, \\
    \\
  \Cset(f)\ \subset\ \Cset\left(\displaystyle\frac{-i(t+i)(1+t)t}{(-t+i)(-1+t)}\right)\ \subset\ \Cset(x)\,.
\end{array}\]

As there is a bijection in this case, the corresponding two decompositions are complete.
\end{exmp}

In order to compute the fixing group of a function $f$ we can solve the system of polynomial
equations obtained from
\[f\left({\frac{ax+b}{cx+d}}\right)=f(x).\]
This can be reduced to solving two simpler systems, those given by
\[f(ax+b)=f(x)\qquad\mbox{and}\qquad f\left({\frac{ax+b}{x+d}}\right)=f(x).\]
This method is simple but inefficient; we will describe another method that is faster in practice.

We need to assume that $\K$ has sufficiently many elements. If not, we take an extension of $\K$
and later we check which of the computed elements are in $\Gamma(\K)$ by solving simple systems of
linear equations.

\begin{thm}
Let $f\in\K(x)$ of degree $m$ in normal form and $u=\displaystyle\frac{ax+b}{cx+d}$ such that
$f\circ u=f$.
\item[(i)] $a\neq 0$ and $d\neq 0$.

\item[(ii)] $f_N(b/d)=0$.

\item[(iii)] If $c=0$ (that is, we take $u=ax+b$), then $f_N(b)=0$ and $a^m=1$.

\item[(iv)] If $c\neq 0$ then $f_D(a/c)=0$.
\end{thm}

\begin{pf}

(i) Suppose $a=0$. We can assume $u=1/(cx+d)=(1/x)\circ(cx+d)$. But if we consider $f(1/x)$, its
numerator has smaller degree than its denominator. As composing on the right with $cx+d$ does not
change those degrees, it is impossible that $f\circ u=f$. Also, as the inverse of $u$ is
$\displaystyle\frac{dx-b}{-cx+a}$, we have $d\neq 0$.

(ii) Let
\[f=\frac{a_mx^m+\cdots+a_1x}{b_{m-1}x^{m-1}+\cdots+b_0}.\]
The constant term of the numerator of $f\circ u$ is
\[a_mb^m+a_{m-1}b^{m-1}d+\cdots+a_1bd^{m-1}=d^mf_N(b/d).\]
As $d\neq 0$ by (i), we have that $f_N(b/d)=0$. Alternatively, $0=f(0)=(f\circ
u)(0)=f(u(0))=f(b/d)$.

(iii), (iv) They are similar to the previous item.
\end{pf}

We can use this theorem to compute the polynomial and rational elements of $G(f)$ separately.

\begin{alg}\label{alg-fixing-group} $ $
\begin{description}
    \item \emph{INPUT}: $f\in\K(x)$.

    \item \emph{OUTPUT}: $G(f)=\{w\in\K(x): f\circ w=f\}$.
\end{description}
\begin{description}
    \item \emph{A}. Compute units $u,v$ such that $\overline{f}=u\circ f\circ v$ is in normal
form. Let $m=\deg\,f$. Let $L$ be an empty list.

    \item \emph{B}. Compute $A=\{\alpha\in\K: \alpha^m=1\}$,
$B=\{\beta\in\K: \overline{f}_N(\beta)=0\}$ and $C=\{\gamma\in\K:\ \overline{f}_D(\gamma)=0\}$.

    \item \emph{C}. For each $(\alpha,\beta)\in A\times B$, check if
$\overline{f}(\alpha x+\beta )=\overline{f}(x)$. In that case add $ax+b$ to $L$.

    \item \emph{D}. For each $(\beta,\gamma)\in B\times C$, let
$w=\displaystyle{\frac{c\gamma\,x+\beta}{c\,x+1}}$. Compute all values of $c$ for which
$\overline{f}\circ w=\overline{f}$. For each solution, add the corresponding unit to $L$.

    \item \emph{E}. Let $L=\{w_1,\ldots,w_k\}$. Return $\{v\circ w_i\circ v^{-1}: i=1,\ldots,k\}$.
\end{description}
\end{alg}

\textbf{Analysis}. It is clear that the cost of the algorithm heavily depends on the complexity of
the best algorithm to compute the roots of a univariate polynomial in the given field. We analyze
the bit complexity when the ground field is the rational number $\Qset$. We will use several
well--known results about complexity, those can be consulted in the book \cite{GG99}.

In the following, $M$ denotes a multiplication time, so that the product of two polynomials in
$\K[x]$ with degree at most $m$ can be computed with at most $M(m)$ arithmetic operations. If $\K$
supports the Fast Fourier Transform, several known algorithms require $O(n\log n\log\log n)$
arithmetic operations. We denote by $l(f)$ the maximum norm of $f$, that is, $l(f)= \|f\|_\infty =
\max |a_i| $ of a polynomial $f =\sum_i a_ix^i\in \Zset[x]$.

Polynomials in $f, g \in \Zset[x]$ of degree less than $m$ can be multiplied using $O(M(m(l+\log
m)))$ bit operations, where $l=\log \max (l(f), l(g))$.

Now, suppose that the given polynomial $f$ is squarefree primitive, then we can compute all its
rational roots with an expected number of $T(m,\log l(f))$ bit operations, where $T(m,\log l(f))=$

\[\begin{array}{l}
  O(m\log(ml(f)) (\log^2\log\log m+(\log\log l(f))^2\log\log\log l(f)) \\
  +m^2M(\log(ml(f))).
\end{array}\]

We discuss separately the algorithm steps. Let $f=f_N/f_D$, where $f_N, f_D \in \Zset[x]$ and let
$l=\log \max (l(f_N), l(g_D))$ and $m =\deg f$.

\textbf{Step A.} Let $u\in \Qset(x)$ be a unit such that $g_N/g_D=u(f)$ with $ \deg g_N > \deg
g_D$. Such a unit always exists:
\begin{description}
\item [--] If $\deg f_N= \deg f_D$. Let $u= 1/(x-a)$, where $a \in \Qset$ verifies $\deg f_N -a\deg f_D < \deg f_N$.
\item [--] If $\deg f_N < \deg f_D$, let $u=1/x$.
\end{description}
Now, let $b \in \Zset$ such that $g_D(b) \neq 0$. Then $h_N/h_D=g_N(x+b)/g_D(x+b)$ verifies
$h_D(0)\neq 0$ and the rational function $(x-h(0)) \circ h_N/h_D$ is in normal form. Obviously,
the complexity in this step is dominated on choosing $b$. In the worst case, we have to evaluate
the integers $0,1,\ldots, m$ in $g_D$. Clearly, a complexity bound is $O(M(m^3l))$.

 \textbf{Step B.} Compute the set $A$ can be done on constant time. Now, in order to
 compute the complexity, we can can suppose, without loss of generality, that
 $\overline{f_N}$ and $\overline{f_D}$ are squarefree and primitive.
 Then the bit complexity to compute both set $B$ and set $C$ is
$T(m, ml)$.

\textbf{Step C.} A bound for the cardinal of Ê$A$ is $4$ and $m$ for the cardinal of $B$. Then, we
need to check $4m$ times if $\overline{f}(\alpha x+\beta )=\overline{f}(x)$ for each each
$(\alpha,\beta)\in A\times B$. So, the complexity of this step is bounded by $O(M(m^4l))$.

\textbf{Step D.} In the worst case the cardinal of $B \times C$ is $m^2$. This step requires to
compute all rational roots of $m^2$ polynomials $h(x)$ given by the equation:
\[\overline{f}\circ w=\overline{f},\]
for each $(\beta,\gamma)\in B\times C$, where $w=\displaystyle{\frac{c\gamma\,x+\beta}{c\,x+1}}$.
A bound for the degree of $h(x)$ is $m^2$. The size of the coefficients is bounded by $ml$, so a
bound for total complexity of this step is $m^4 T(m^2,lm^2)$.

\textbf{Step E.} Finally, this step requires substituting at most $2m$ rational functions of
degree $m$ and the coefficients size is bounded by $lm^3$. So, abound for the complexity is
$O(M(m^4l))$.

We can conclude that the complexity of this algorithm is dominated by that of step D, that is,
$m^4T(m^2,lm^2)$. Of course, a worst bound for this is $O(m^8l^2)$.

The following example illustrates the above algorithm:

\begin{exmp}
Let
\[f=\frac{(-3x+1+x^3)^2}{x(-2x-x^2+1+x^3)(-1+x)}\in\Qset(x).\]

We normalize $f$: let $u=\displaystyle{\frac{1}{x-9/2}}$ and $v=\displaystyle{\frac{1}{x}}-1$,
then
\[\overline{f}=u\circ f\circ v=\frac{-4x^6-6x^5+32x^4-34x^3+14x^2-2x}{27x^5-108x^4+141x^3-81x^2+21x-2}\]
is in normal form.

The roots of the numerator and denominator of $\overline{f}$ in $\Qset$ are $\{0,1,1/2\}$ and
$\{1/3,2/3\}$ respectively. The only sixth roots of unity in $\Qset$ are $1$ and $-1$; as
$\mathrm{char}\ \Qset=0$ there cannot be elements of the form $x+b$ in $G(\overline{f})$. Thus,
there are two polynomial candidates: $-x+1/3$, $-x+2/3$. A quick computation reveals that none of
them fixes $\overline{f}$.

Let $w=\displaystyle\frac{c\beta\,x+\alpha}{c\,x+1}$. As $\alpha\in\{0,1,1/2\}$ and
$\beta\in\{1/3,2/3\}$, another quick computation shows that
\[G(\overline{f})=\left\{\ x,\quad \frac{-x+1}{-3x+2},\quad \frac{-2x+1}{-3x+1}\ \right\}\]
and
\[G(f)=v\cdot G(\overline{f})\cdot v^{-1}=\left\{\ x,\quad \frac{1}{1-x},\quad \frac{x-1}{x}\ \right\}.\]

From this group we can compute a proper component of $f$ as in the proof of Theorem
\ref{conds-indesc}, obtaining $f=g(h)$ with

\[h=\frac{-3\,x+1+x^3}{(-1+x)x}, \quad g=\frac{x^2}{x-1}.\]

\end{exmp}

In the next section we will use these tools to investigate the
number of components of a rational function.

\section{Ritt's Theorem and number of components}

One of the classical Ritt's Theorems (see \cite{Rit22}) describes the relation among the different
decomposition chains of a tame polynomial. Essentially, all the decompositions have the same
length and are related in a rather simple way.

\begin{defn}
A \emph{bidecomposition} is a 4-tuple of polynomials $f_1,g_1,f_2,g_2$ such that $f_1\circ
g_1=f_2\circ g_2$, $\deg\,f_1=\deg\,g_2$ and $(\deg\,f_1,\deg\,g_1)=1$.
\end{defn}

\begin{thm}[Ritt's First Theorem]\label{ritt1}
Let $f\in\K[x]$ be tame and
\[f=g_1\circ\cdots\circ g_r=h_1\circ\cdots\circ h_s\]
be two complete decomposition chains of $f$. Then $r=s$, and the sequences
$(\deg\,g_1,\ldots,\deg\,g_r)$, $(\deg\,h_1,\ldots,\deg\,h_s)$ are permutations of each other.
Moreover, there exists a finite chain of complete decompositions
\[f=f_1^{(j)}\circ\cdots\circ f_r^{(j)}\,,\quad j\in\{1,\ldots,k\}\,,\]
such that
\[f_i^{(1)}=g_i\,,\quad f_i^{(k)}=h_i\,,\quad i=1,\ldots,r\,,\]
and for each $j<k$, there exists $i_j$ such that the $j$-th and
$(j+1)$-th decomposition differ only in one of these aspects:

\item[(i)] $f_{i_j}^{(j)}\circ f_{i_j+1}^{(j)}$ and $f_{i_j}^{(j+1)}\circ f_{i_j+1}^{(j+1)}$ are
equivalent.

\item[(ii)] $f_{i_j}^{(j)}\circ f_{i_j+1}^{(j)}=f_{i_j}^{(j+1)}\circ f_{i_j+1}^{(j+1)}$ is a bidecomposition.
\end{thm}

\begin{pf}
See \cite{Rit22} for $\K=\Cset$, \cite{Eng41} for characteristic zero fields and \cite{FM69b}, \cite{Sch00} for
the general case. \qed
\end{pf}

Unlike for polynomials, it is not true that all complete decompositions of a rational function
have the same length, as shown in Example \ref{ej-fixed-field2}. The paper \cite{GS05} presents a
detailed study of this problem for non tame polynomial with coefficients over a finite field. The
problem for rational functions is strongly related to the open problem of the classes of rational
functions which commute with respect to composition, see \cite{Rit23}. In this section we will
give some ideas about the relation between complete decompositions and subgroup chains that appear
by means of Galois Theory.

Now we present another degree 12 function, this time with coefficients in $\Qset$, that has two
complete decomposition chains of different length. This function arises in the context of
Monstrous Moonshine as a rational relationship between two modular functions (see for example the
classical \cite{CN79} for an overview of this broad topic, or the reference \cite{McS04}, in
Spanish, for the computations in which this function appears).

\begin{exmp}
Let $f\in\Qset(x)$ be the following degree 12 function:
\[f=\frac{x^3(x+6)^3(x^2-6\,x+36)^3}{(x-3)^3(x^2+3\,x+9)^3}.\]

$f$ has two decompositions:
\[\begin{array}{rl}
  f\ =  &  g_1\circ g_2\circ g_3\ =\ x^3\ \circ\ \displaystyle\frac{x(x-12)}{x-3}\ \circ\ \displaystyle\frac{x(x+6)}{x-3}\ =  \\
    \\
  =  &  h_1\circ h_2\ =\ \displaystyle\frac{x^3(x+24)}{x-3}\ \circ\ \displaystyle\frac{x(x^2-6\,x+36)}{x^2+3\,x+9}.
\end{array}\]

All the components except one have prime degree, hence are indecomposable; the component of degree
4 cannot be written as composition of two components of degree 2.
\end{exmp}

If we compute the groups for the components on the right in $\Qset$ we have:
\[\begin{array}{l}
  G_\Qset(f)\ =\ G_\Qset(g_2\circ g_3)\ =\ G_\Qset(g_3)\ =\ \left\{\,\displaystyle\frac{3x+18}{x-3},\ x\,\right\}, \\
  G_\Qset(h_2)\ =\ \{x\}.
\end{array}\]
However, in $\Cset$:
\[\begin{array}{rcl}
  G_\Cset(f)  &  =  &  \left\{\ \displaystyle{\frac{3\alpha_ix+18\alpha_i}{x-3}\,,\quad \frac{3\alpha_ix-18-18\alpha_i}{x-3\alpha_i}\,,\quad \frac{3\alpha_ix+18}{x+3\alpha_i+3}\,,}\right.  \vspace{1ex}  \\
    &  &  \quad \left.\displaystyle\frac{3x+18\alpha_i}{x-3\alpha_i}\,,\quad \displaystyle\frac{3x+18}{x-3}\,,\quad \alpha_ix\,,\ x\ \right\},  \\
    \\
  G_\Cset(g_2\circ g_3)  &  =  &  \left\{\ \displaystyle\frac{3\alpha_ix-18-18\alpha_i}{x-3\alpha_i}\,,\quad \frac{3x+18}{x-3}\,,\ x\ \right\},  \vspace{1ex}  \\
  G_\Cset(g_3)  &  =  &  \left\{\ \displaystyle\frac{3x+18}{x-3}\,,\ x\ \right\},  \vspace{1ex}  \\
  G_\Cset(h_2)  &  =  &  \left\{\ \displaystyle\frac{3\alpha_ix+18}{x+3\alpha_i+3}\,,\ x\ \right\}
\end{array}\]
where $\alpha_i,\ i=1,2$ are the two non-trivial cubic roots of
unity.

In order to obtain the function in Example \ref{ej-fixed-field2}, we used Theorem \ref{props-fix},
and in particular the existence of a bijection between the subgroups of $A_4$ and the intermediate
fields of a function that generates the corresponding field. The existence of functions with this
property has been known for some time, as its construction from any group isomorphic to $A_4$ is
straightforward. On the other hand, the example above is in $\Qset(x)$, but there is no bijection
between groups and intermediate fields.

In general, there are two main obstructions for this approach. On one hand, there is no bijection
between groups and fields in general, as the previous example shows for $\Qset$. On the other
hand, only some finite groups can be subgroups of $\mathrm{PGL}_2(\K)$. The only finite subgroups
of $\mathrm{PGL}_2(\Cset)$ are $C_n$, $D_n$, $A_4$, $S_4$ and $A_5$, see \cite{Kle56}. In fact,
this is true for any algebraically closed field of characteristic zero (it suffices that it
contains all roots of unity). Among these groups, only $A_4$ has subgroup chains of different
length. This is even worse if we consider smaller fields as the next known result shows:

\begin{thm}
Every finite subgroup of $\mathrm{PGL}_2(\Qset)$ is isomorphic to either $C_n$ or $D_n$ for some
$n\in\{2,3,4,6\}$.
\end{thm}

Indeed these all occur, unfortunately none of them has two subgroup chains of different lengths,
so no new functions can be found in this way.

\section{Conclusions}

In this paper we have presented several counterexamples to the generalization of the first Ritt
theorem to rational functions. We also introduced and analyzed several concepts of Galois Theory
that we expect to be interesting in providing more structural information in this topic. Also, we
show a use of techniques from Computational Algebra results to find bounds for the size of a field
that contains all decompositions of a given function; we expect that general properties of
Gr\"obner bases can be applied to this end in order to obtain general bounds.

\end{document}